\documentstyle[art11]{article}

\newdimen\mypt
\newdimen\halfpt
\def\pic#1 #2 #3;{\raisebox{0.6ex}{\raisebox{-#2\halfpt}{\vbox
      to #2\mypt{\hbox to #3\mypt{\special{em:graph
      #1.bmp}\hfill}\vfill}}}}

\def\proved{\ifmmode\eqno\Box\medskip\else
        \nobreak\hfill\nopagebreak\discretionary{}
        {\hbox to\textwidth{\hfill$\Box$}}{\hbox{$\Box$}}\par
        \addvspace\medskipamount\fi}

\newcounter{defin}

\newcounter{remark}

\font\msamfive=msam5 at 3pt
\def\uparc(#1_#2#3_#4){{\mathop{#1_{#2}#3_{#4}}\limits^{\displaystyle
\frown\kern-3pt\llap{\raisebox{1.5pt}{\msamfive\char121}}
\phantom{\scriptscriptstyle#4}
}}}
\def\loarc(#1_#2#3_#4){{\mathop{#1_{#2}#3_{#4}}\limits_{\displaystyle
\smile\kern-3pt\llap{\raisebox{1.7pt}{\msamfive\char113}}
\phantom{\scriptscriptstyle#4}
}}}

\def\sing{{\rm sing \,}}

\author{S. D. Tyurina\thanks{This work is partly
supported by INTAS, grant no.~97-1644.}
}
\date{Dep. of Phys. and Math.\\
Kolomensky Pedagogical Institute\\
Kolomna 140411 RUSSIA\\
e-mail: tyurina@mccme.ru}
\title{{\sc On formulas of the Lannes and Viro-Polyak type for finit degree
invariants.}\\
{\scriptsize(It is very brief text. More complete version of this paper
is in Mat.Zametki {\bf t.65} N6 (1998).
)}}

\begin{document}
\baselineskip=32pt
\maketitle
UDK 515.162
\abstract{
Vassiliev's knot invariants can be computed in different ways but many of them
as Kontsevich integral are very difficult. We consider more visual diagram
formulas of the type Polyak-Viro and give new diagram formula for the two
basic Vassiliev invariant of degree 4.}
\section{The Vassiliev knot invariants.}
Let $K:S^1 \longrightarrow {\bf R^3}$
be an oriented knot and $K^{\sing}_n:S^1
\longrightarrow {\bf R^3}$ be a singular knot
with $n$ double points.
Denote by ${\cal K}$ the space of knots and
by ${\cal K}^{\sing}$ the space of singular knots.
Any knot invariant $V:{\cal K} \longrightarrow {\bf Q}$ may be extended from
ordinary knots to singular knots by next inductive rule:

\begin{picture}(420,22)
\put(80,9){$V^{(i)}($}
\put(120,11){\circle*{2}}
\put(132,9){$)=V^{(i-1)}($}
\put(217,9){$)-V^{(i-1)}($}
\put(299,9){).}
\put(420,9){(1)}
\put(109,0){\vector(1,1){22}}
\put(131,0){\vector(-1,1){22}}
\put(193,0){\vector(1,1){22}}
\put(215,0){\line(-1,1){10}}
\put(203,12){\vector(-1,1){10}}
\put(274,0){\line(1,1){10}}
\put(288,12){\vector(1,1){10}}
\put(297,0){\vector(-1,1){22}}
\end{picture}

{\bf Definition.} A knot invariant $V:{\cal K} \longrightarrow
{\bf Q}$ is said to be Vassiliev invariant of degree
less than or equal to $n$, if there exists $n \in {\bf N}$ such that
$$V^{(n+1)} \equiv 0.$$

\section{Weight systems.}
For any singular knot we construct it's chord diagram.
A chord diagram $D_n$ of the singular knot $K^\sing$ is the circle with
pre-images of double points connected with  chords.

Denote by ${\cal D}_n$ the free abelian group generated by diagrams with
$n$ chord. For any invariant $f$ with values in an abelian group
$A$ there exists a homonorphism $W_f:{\cal D}_n \longrightarrow A,$ given by
$W_f = f|_{{\cal K}^{sing}_n},$ where ${\cal K}^\sing_n$ is the space of
knots with $n$ double points.

{\bf Definition.} A linear function $W$ is called a weight system of degree $n$
if it satisfies next relations:

\begin{picture}(420,20)
\put(0,9){1-term relation:}
\put(100,9){$W_n($}
\put(135,11){\circle{15}}
\put(147,9){$)=0$}
\put(135,3){\line(0,1){15}}
\end{picture}

\noindent (this diagram has $n$ chords, one of which is isolated) and

\begin{picture}(420,20)
\put(0,9){4-term relation:}
\put(100,9){$W_n($}
\put(135,11){\circle{15}}
\put(148,9){$)-W_n($}
\put(200,11){\circle{15}}
\put(213,9){$)+W_n($}

\put(265,11){\circle{15}}
\put(280,9){$)-W_n($}

\put(332,11){\circle{15}}
\put(347,9){$)=0$}

\put(134,4){\line(-1,2){5}}
\put(136,4){\line(1,2){5}}

\put(197,4){\line(1,1){10}}
\put(203,4){\line(-1,1){10}}
\put(263,4){\line(1,2){6.5}}
\put(271,14){\line(-1,0){12}}
\put(329,4){\line(1,1){10}}
\put(338,16){\line(-1,0){11.5}}
\end{picture}

\noindent these diagrams have $n$ chords, $(n-2)$ of  which are not drown
here and 2 chords are positioned as shown.

Defineded above homonorphism $W_f$ is a weight system.

Examples 1. Weight systems of degrees 2 and 3.

\noindent
\begin{picture}(400,40)
\put(0,22){$W_2($}
\put(0,2){$W_2($}
\put(30,25){\circle{15}}
\put(30,7){\circle{15}}
\put(40,22){)=1}
\put(40,2){)=0}
\put(25,20){\line(1,1){10}}
\put(35,20){\line(-1,1){10}}
\put(25,1){\line(0,1){11}}
\put(35,1){\line(0,1){11}}
\put(200,22){$W_3($}
\put(200,2){$W_3($}
\put(230,25){\circle{15}}
\put(230,7){\circle{15}}
\put(240,22){)=2}
\put(240,2){)=1}
\put(225,20){\line(1,1){10}}
\put(235,20){\line(-1,1){10}}
\put(225,1){\line(0,1){11}}
\put(235,1){\line(0,1){11}}
\put(222,25){\line(1,0){15}}
\put(222,7){\line(1,0){15}}
\end{picture}

\hskip6cm $W_3=0$ in other cases.

\section{"Coordinates" on knots.}
Let $D$ be a diagram of knot $K$ and $x$ be a double point of $D.$
Numerate branches in neighbourhood of $x$ according to the order of their
passing. Define function $\delta_x$ by next rule:

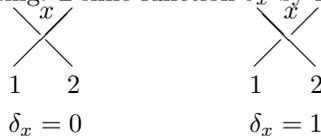
\begin{picture}(400,32)
\put(143,10){\line(1,1){22}}
\put(165,10){\line(-1,1){10}}
\put(141,0){$1$}
\put(163,0){$2$}
\put(141,-15){$\delta_x = 0$}
\put(153,22){\line(-1,1){10}}
\put(152.5,28){$x$}
\put(234,10){\line(1,1){10}}
\put(248,22){\line(1,1){10}}
\put(245,28){$x$}
\put(257,10){\line(-1,1){22}}
\put(232,-15){$\delta_x = 1$}
\put(255,0){$2$}
\put(232,0){$1$}

\end{picture}

\vskip0.5cm
(the orientations of branches are not important).

Define function $\varepsilon_x$ as follows:

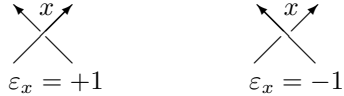
\begin{picture}(400,32)
\put(143,10){\vector(1,1){22}} 
\put(165,10){\line(-1,1){10}}
\put(153,22){\vector(-1,1){10}}
\put(234,10){\line(1,1){10}}
\put(248,22){\vector(1,1){10}}
\put(257,10){\vector(-1,1){22}}
\put(141,0){$\varepsilon_x = +1$}
\put(232,0){$\varepsilon_x = -1$}
\put(152.5,28){$x$}
\put(245,28){$x$}
\end{picture}

Examples 2. "Coordinates" on trefoil.

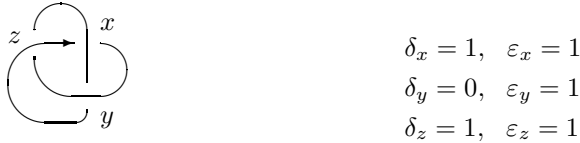
\begin{picture}(400,55)
\put(20,35){\vector(1,0){5}}
\put(20,40){\oval(20,20)[t]}
\put(30,20){\line(0,1){20}}
\put(25,20){\oval(50,30)[tl]}
\put(30,30){\oval(40,30)[bl]}
\put(30,15){\line(1,0){5}}
\put(15,10){\oval(30,10)[br]}
\put(35,25){\oval(20,20)[r]}
\put(15,20){\oval(30,30)[bl]}
\put(35,40){$x$}
\put(35,5){$y$}
\put(0,35){$z$}
\put(150,30){$\delta_x=1,\,\,\,\, \varepsilon_x=1$}
\put(150,15){$\delta_y=0,\,\,\,\, \varepsilon_y=1$}
\put(150,0){$\delta_z=1,\,\,\,\, \varepsilon_z=1$}
\end{picture}

\section{Formulas of Lannes.}
For the basic Vassiliev invariants of degree 2 and 3, which take values 0
on trivial knot and 1 on trefoil we have next formulas:
$$V_2(K)=1/2 \sum_{\{x,y\}\in P_2}(-1)^{\delta_x + \delta_y} W_2(\{x,y\})
\varepsilon_x \varepsilon_y [\delta_x (1- \delta_y )+\delta_y (1-\delta_x)],$$
$$V_3(K)=1/2 \sum_{\{x,y,z\}\in P_3} (-1)^{\delta_x+\delta_y+\delta_z}
W_3(\{x,y,z\}) \varepsilon_x \varepsilon_y \varepsilon_z
[\delta_y (1- \delta_x)(1- \delta_z)- \delta_x \delta_z (1- \delta_y)],$$
where the sum is taken over all unordered pairs (triplets) of double points
of planar projection, $W_2(\{x,y\})$ ($W_3(\{x,y,z\})$) is weight of chord
diagram corresponding to pair (triplet) of double points.

\section{The Gauss diagrams.}
A chord diagram of the singular knot is the circle with pre-images
of double points connected with  chords.
To obtain the analogus diagram
of an ordinary knot (that is called an arrow diagram) from the chord diagram
of corresponding singular knot we must give the information on overpasses and
underpasses. Each chord is oriented from the upper branch to the lower one
and equipped with the sign (the local writh number of corresponding double
point of planar projection of the knot).

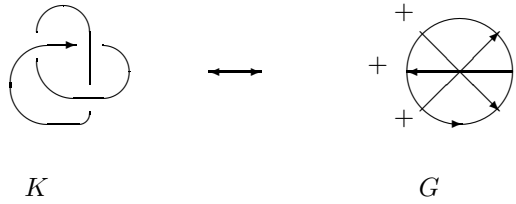
\begin{picture}(400,65)
\put(100,45){\vector(1,0){5}}
\put(100,50){\oval(20,20)[t]}
\put(110,30){\line(0,1){20}}
\put(105,30){\oval(50,30)[tl]}
\put(110,40){\oval(40,30)[bl]}
\put(110,25){\line(1,0){5}}
\put(95,20){\oval(30,10)[br]}
\put(115,35){\oval(20,20)[r]}
\put(95,30){\oval(30,30)[bl]}

\put(165,35){\vector(-1,0){10}}
\put(165,35){\vector(1,0){10}}
\put(250,35){\circle{40}}
\put(250,15){\vector(1,0){1}}
\put(270,35){\vector(-1,0){40}}
\put(235,20){\vector(1,1){30}}
\put(235,50){\vector(1,-1){30}}
\put(215,35){$+$}
\put(225,15){$+$}
\put(225,55){$+$}
\end{picture}

\hskip3cm {$K$\qquad \qquad \qquad \qquad \qquad \qquad \qquad $G$}

\section{Formulas of Viro-Polyak.}
Denote by $<A,G>$ algebraic number of subdiagrams of given combinatorial
type $A \quad A \subset G$ and let $<\sum_{i}n_iA_i,G> = \sum_{i}n_i<A_i,G>,
\quad n_i \in {\bf Q}$ by definition.
Then

\begin{picture}(400,22)
\put(193,5){\vector(1,1){11}}
\put(203,5){\vector(-1,1){11}}
\put(198,10){\circle{15}}
\put(198,17){\circle*{2}}
\put(135,5){\it $V_2(K)=<$}
\put(210,5){$,G>,$}
\end{picture}

\begin{picture}(400,22)
\put(183,5){\vector(1,1){11}}
\put(193,5){\vector(-1,1){11}}
\put(188,17){\vector(0,-1){15}}
\put(188,10){\circle{15}}
\put(115,5){\it $V_3(K)=<[$}
\put(260,5){$],G>.$}
\put(205,5){$]+\frac{1}{2}[$}
\put(243,13){\vector(1,0){11}}
\put(253,5){\vector(-1,0){11}}
\put(248,2){\vector(0,1){15}}
\put(248,10){\circle{15}}
\end{picture}

{\bf Theorem.} {\it Let $v_3$ be basic Vassiliev invariant of degree 3, which
take values 0 on trivial knot and 1 on trefoil we have next formula:}

\begin{picture}(400,22)
\put(110,5){\it $v_3(K)=<$}
\put(168,5){\vector(1,1){11}}
\put(178,5){\vector(-1,1){11}}
\put(173,17){\vector(0,-1){15}}
\put(173,10){\circle{15}}
\put(185,5){$+$}
\put(208,15){\vector(-1,-1){11}}
\put(198,15){\vector(1,-1){11}}
\put(203,2){\vector(0,1){15}}
\put(203,10){\circle{15}}
\put(215,5){$+$}
\put(303,15){\circle*{2}}
\put(268,17){\circle*{2}}
\put(233,15){\circle*{2}}
\put(170,16){\circle*{2}}
\put(200,16){\circle*{2}}
\put(310,5){$,G>.$}
\put(231,10){\vector(1,0){15}}
\put(241,4){\vector(0,1){13}}
\put(236,17){\vector(0,-1){14}}
\put(238,10){\circle{15}}
\put(250,5){$+$}
\put(276,10){\vector(-1,0){15}}
\put(272,4){\vector(0,1){13}}
\put(265,17){\vector(0,-1){14}}
\put(268,10){\circle{15}}
\put(280,5){$+$}
\put(306,10){\vector(-1,0){15}}
\put(300,4){\vector(0,1){13}}
\put(295,17){\vector(0,-1){14}}
\put(298,10){\circle{15}}
\end{picture}

Proof see in [5].

{\bf The Vassiliev module.}
The Vassiliev module of degree $n$ is a module over ${\bf Q}$
generated by isotopies classes of oriented knots and singular knots with
next relations:
\newline
1. E=0, where E is the trivial knot, \newline
2. the Vassiliev skein-relation (1), \newline
3. $K_m^{sing}=0$, if $m>n.$

{\bf Theorem.} Any knot $K$ in the Vassiliev module of degree $n \quad
(n>1)$ has following expansion:
$$K=\sum_{i=1}^{r+s}v_i(K)K_i,$$
where $r$ is the dimension of the space of the Vassiliev invariants of degree
less than or equal to $(n-1),$ \newline
$s$ is the dimension of the space of weight system of degree $n,$ \newline
$v_i$ is the Vassiliev invariants of degree less than or equal to $n,$ \newline
$K_i$ is fixed basic knots.
\bigskip
Examples. \newline
$n=2$ (J.Lannes)
$$K=V_2(K)T$$
$n=3$ (J.Lannes)
$$K=V_3(K)[T + H] - V_2(K)H$$
\newline $n=4$
$$K=[V_2/2+V_3/2-3V_4^1+4V_4^2+V_4^3](K)T+[V_2/2-V_3/2+V_4^3](K)T^*+$$
$$[-V_4^1+V_4^2+2V_4^3](K)H+[V_4^1-3V_4^2](K)F+V_4^2(K)P$$


\newpage
\begin{thebibliography}{99}
\bibitem{Lannes}
J.Lannes.
Sur les invariants de Vassiliev de degre$\acute e$ inferieur
ou $\acute e$gal $\grave a$ 3.
{\it L'En\-seignement Math$\acute e$matique} {\bf t.39} (1993), pp.~295--316.
\bibitem{VP}
M.Polyak, O.Viro.
Gauss diagram formulae for Vassiliev invariants.
{\it Int. Math. Res. Notices} {\bf 11} (1994), pp.~445-453.
\bibitem{Vas}
V.A.Vassiliev.
Complemets of discriminants of smoos maps.
{\it Amer. Math. Soc. Transl.} (1998).
\bibitem{Kon}
M.Kontsevich.
Vassiliev's knot invariants.
{\it Adv.in Sov.Math.} {\bf 16} (1993), pp.137-150.
\bibitem{Tyur}
S.Tyurina.
On formulas of the Lannes and Viro-Polyak type for finit degree invariants.
{\it Mat. zametki} {\bf 65} N 6 (1999).
\end{thebibliography}
\end{document}